\documentclass[11pt]{amsart}

\title{Modular forms of weight one: Galois representations and dimension.\\
\tiny{Notes des expos\'es 10, 12 et 13 du s\'eminaire Frank Thorne.}}
\author{Denis Trotabas}
\date{}
\usepackage{amsmath,latexsym,amssymb,amsfonts,mathrsfs,bbm,fancyhdr}
\usepackage[colorlinks=true, linkcolor=black,citecolor=black,
urlcolor=black]{hyperref}
\newcommand{\e}{{\rm e}}
\newcommand{\Q}{\mathbf{Q}}
\newcommand{\adele}{\mathbbmss{A}}

\newcommand{\HH}{\mathfrak{H}}
\newcommand{\N}{\mathbf{N}}
\newcommand{\Z}{\mathbf{Z}}
\newcommand{\F}{\mathbf{F}}
\newcommand{\R}{\mathbf{R}}
\newcommand{\C}{\mathbf{C}}
\newcommand{\entier}{\mathcal{O}_K}
\newcommand{\n}{\mathfrak{n}}
\newcommand{\p}{\mathfrak{p}}

\newcommand{\Gal}[1] {{ Ga\ell}(\overline{#1}/#1) } 
\newcommand{\GL}{{\rm GL}}

\def\build#1_#2^#3{\mathrel{\mathop{\kern 0pt#1}\limits_{#2}^{#3}}} %CLENET
\def\equi_#1{\build{\sim}_#1^{}}
\def\Equi#1#2{\equi_{\scriptscriptstyle#1\to #2}}

\newtheorem{probleme}{\sc Problem}[section]
\newtheorem{theoreme}{\sc Theorem}[section]

\newtheorem{conj}{\sc{Conjecture}}[section]
\newtheorem{lemme}{\sc{Lemma}}[section]
\newtheorem{proposition}{\sc{Proposition}}[section]
\newtheorem{corollaire}{\sc{Corollary}}[theoreme]

\begin{document}

\begin{abstract}
The present notes are the expanded and polished version of three lectures given in Stanford, concerning the analytic and arithmetic properties of weight one modular forms. The author tried to write them in a style accessible to non-analytically oriented number theoritists: in particular, some effort is made to be precise on statements involving uniformity in the parameters. On the other hand, another purpose was to provide an introduction, together with a set of references, consciously kept small, to the realm of Galois representations, for non-algebraists -- like the author. The proofs are sketched, at best, but we tried to motivate the results, and to relate them to interesting conjectures. The author thanks warmly Daniel Bump, Brian Conrad, Kannan Soundararajan, Akshay Venkatesh and Justin Walker for helpful comments and inspiring conversations.
\end{abstract}

\maketitle

\section{Some notations.}
Let $\HH=\{z\in\C; \Im(z)>0\}$ be the Poincar\'e upper-half plane. 

Let $k$ and $q$ be two integers, and, as usual, $\Gamma_0(q)$ be the subgroup of ${\rm SL}(2,\Z)$ of matrices whose lower left entries are divisible by $q$. It acts on $\HH$ by fractional linear transformations: $\bigl[\begin{smallmatrix} a&b\\c&d\end{smallmatrix}\bigr] \cdot z=\frac{az+b}{cz+d}\,$. 

Let $\chi$ be a Dirichlet character modulo $q$: it defines a character on $\Gamma_0(q)$, by evaluating $\chi$  at the lower right entry.\\

If $X$ is a finite set, $|X|$ denotes its cardinality; we reserve the letters $p,\ell$ for prime numbers, and $n,m$ for integers; $\pi(x)$ is the number of primes less than $x$. Recall that $\pi(x)\Equi{x}{\infty}\frac{x}{\log(x)}\,${}.

The letters $K,E,k$ (resp. $K_\lambda$) denote fields (resp. the completion of $K$ with respect to the valuation associated to $\lambda$), and $\entier, \mathcal{O}_\lambda$ stand for the rings of integers of $K, K_\lambda$ in the relevant situations.\\

For a complex number $z$, the notation $\e(z)$ stands for $\exp(2\pi i z)$.\\

The notation $f(x)\ll_A g(x)$ means that for any $A$, there exists a real number $C(A)$ such that for any $x$, $|f(x)|\leq C(A)\cdot |g(x)|$; if one adds ``as $x\to\infty$'', it means that the last inequality holds for $x\geq x(A)$ for some real number $x(A)$. In the same spirit, the notation $f(x)=o_{x\to x_0}(g(x))$ (resp. $f(x)=\mathcal{O}_{x\to x_0}(g(x))$) means that the quotient $f(x)/g(x)$ is defined in a (pointed) neighbourhood of $x_0$, and that $|f(x)/g(x)|$ tends to zero (resp. stays bounded) when $x$ tends to $x_0$. 

Sometimes, we used the notation $x\doteq y$ to mean $x=cy$ for some constant $c$, specifying the parameters on which $c$ may depend; in these instances, the reader can check that the corresponding identities are not as trivial as this notation suggests.
\section{Modular forms}

\subsection{} For any holomorphic function $f$ defined on $\HH$ and $\gamma\in\Gamma_0(q)$, we define:
$$f_{|_\gamma}(z)=\chi(\gamma)^{-1}(cz+d)^{-k}f(\gamma(z))$$
Consider the following properties:\\
(M1): For any $\gamma\in\Gamma_0(q)$, $f_{|_\gamma}=f$.\\
This implies, by Fourier analysis, that for any $\sigma\in{\rm SL}(2,\Z)$ one has a decomposition:
$$f_{|_\gamma}(z)=\sum_{n\in\Z}c_n(f,\sigma)\exp(2i\pi nz)$$
The holomorphy at $i\infty$ is then expressed by:\\
(M2): For any $\sigma\in{\rm SL}(2,\Z)$, $c_n(f,\sigma)=0$ for all negative $n$.\\
(M2'): For any $\sigma\in{\rm SL}(2,\Z)$, $c_n(f,\sigma)=0$ for all non-positive $n$.\\

\subsection{} The space of modular forms of weight $k$, level $q$ and nebentypus $\chi$ is the set of holomorphic functions satisfying (M1) and (M2) above; the subspace of modular forms satisfying (M2') as well is called the space of cusp forms, noted $\mathcal{S}_k(q,\chi)$. The latter is finite dimensional (as is the first), and equipped with the Petersson inner product, invariant under the group action (it is a quotient of a Haar measure):
$$\left<f,g\right>=\int_{\Gamma_0(q)\backslash \HH}f(x+iy)\overline{g(x+iy)}y^k\frac{dxdy}{y^2}$$
Note right now that by taking $\gamma=-{\rm I}$, (M1) gives $f(z)=(-1)^k\chi(-1)f(z)$, so if $\chi$ and $k$ don't have the same parity, the space of modular forms is $\{0\}$; we shall exclude this case.\\

\subsection{Hecke operators} On the space of modular forms of weight $k$ and level $q$, one has the so-called Hecke operators, defined as follows for any prime $p$:
\begin{itemize}
\item If $(p,q)=1$, $T_p(f)(z)=\sum_nc_{pn}(f)\e(nz)+\chi(p)p^{k-1}\sum_n c_n(f)\e(pnz)$ -- $p$ is called a good prime.
\item If $p|q$, $T_p(f)(z)=\sum_nc_{pn}(f)\e(nz)$ -- $p$ is a bad prime.
\end{itemize}
The Hecke operators preserve the space of cusp forms; the Hecke operators at good primes all commute, and are normal with respect to the Petersson inner product. These important facts are explained in Miyake \cite{M1}, as are the multiplicativity relations. In particular, if $f$ is an eigenfunction for all the Hecke operators at good primes,  with eigenvalues $\{a_p(f)\}$, one has $c_p(f)=c_f(1)a_p(f)$ at good $p$. To diagonalize further the Hecke operators, and get a good definition of $L$-series, it is necessary to introduce

\subsection{Newforms and oldforms.} Suppose $\chi$ defines a Dirichlet character modulo $q'$, for $q'|q$. For any cusp form $g$ in $\mathcal{S}_k(q',\chi)$, one checks easily that $z\mapsto g(dz)$ defines an element of $\mathcal{S}_k(q,\chi)$, for any $d|q/q'$. Let
$$\mathcal{S}_k^{old}(q,\chi)=\bigcup _{\substack{\chi {\textrm{\tiny factors through }}q'|q\\d|q/q'}}\{z\mapsto g(dz)\,\,:\,\,g\in\mathcal{S}_k(q',\chi)\}$$
be the space of oldforms, and let 
$$\mathcal{S}_k^{new}(q,\chi)=\mathcal{S}_k^{old}(q,\chi)^\perp$$
be the space of newforms (it may be zero!). Then it can be shown that the whole Hecke algebra (i.e. including bad primes) can be diagonalized on the space of newforms. The  primitive Hecke eigenforms (those with $c_1(f)=1$) have distinct eigenvalues outside a finite number of primes (``multiplicity one'', well explained in the adelic setting by Casselman \cite{C}, cf. Gelbart \cite{G} as well). Their $L$-series have an Euler product, absolutely convergent if $\Re(s)>1+k/2$:
$$L(s,f):=\sum_n \frac{a_n(f)}{n^{s}}=\prod_p L(s,f_p)$$
with
\begin{align*}
L(s,f_p)&=\bigg(1-a_p(f)p^{-s}+\chi(p)p^{k-1-2s}\bigg)^{-1}\\&=\bigg(1-\alpha_1(p,f)p^{-s}\bigg)^{-1}\bigg(1-\alpha_2(p,f)p^{-s}\bigg)^{-1}
\end{align*}
at a good prime $p$, and
$$L(s,f_p)=\bigg(1-a_p(f)p^{-s}\bigg)^{-1}$$
at a bad prime, 
along with an analytic continuation (easy to see with the Mellin transform), functional equation -- cf. Bump \cite{Bu}, Miyake \cite{M1}, Iwaniec \cite{I}, etc.

When one proves a theorem, one can often reduce it to the case of new-forms, thanks to this decomposition; be aware it might be tacitly done in what follows...

\subsection{Ramanujan conjecture.} Let $f$ be a primitive newform. The Ramanujan conjecture is the following inequality:
$$|a_p(f)|\leq 2p^{\frac{k-1}{2}}$$ 
which is equivalent to $|\alpha_i(p,f)|= p^{\frac{k-1}{2}}$. It has been a theorem for 35 years now, proven by Deligne for weight greater than two. We will prove below the case $k=1$ (cf. corollary \ref{ramanujan}).

\subsection{Rationality properties.} Let $f\in\mathcal{S}_k(q,\chi)$ be a eigenform for all the Hecke operators at good primes, with Hecke eigenvalues $\{a_f(p)\}_{p\not \,| q}$. Then:
$$\Q(f):=\Q(a_f(p), \chi(p): p\not |q)$$ is a finite extension of $\Q$, and all the Hecke eigenvalues are integers in this extension. If the nebentypus is trivial, then this extension in totally real. Serre explains this cohomologically in his Durham lectures.

\subsection{} An interesting problem is the evaluation of the dimension of the space of cusp forms, when one or more of the parameters ($k,q$) vary. For instance, using Eichler-Selberg trace formula one can prove that (see Knightly-Li \cite{KL} theorem 29.5):
\begin{eqnarray}\label{1}
{\dim}(\mathcal{S}_k(q,\chi))=\frac{k-1}{12}\psi(q)+\mathcal{O}\bigg(q^{1/2}\tau(q)\bigg)
\end{eqnarray}
uniform in $k\geq 2$ and $q$, where $\psi(q)=q\prod_{p|q}(1+p^{-1})$ and $\tau(q)$ is the number of divisors of $q$.\\ 

Similarly, one can bound the dimension of the space of new-forms using the Petersson trace formula (Iwaniec-Luo-Sarnak \cite{ILS}), and one has a uniform estimate for $q$ squarefree, $k\geq 2$:

\begin{eqnarray}
{\dim}(\mathcal{S}^{new}_k(q,\chi))=\frac{k-1}{12}\varphi(q)+\mathcal{O}\bigg( (kq)^{2/3}\bigg)
\end{eqnarray}
with $\varphi(q)=q\prod_{p|q}(1-p^{-1})$ the Euler phi ([fi]) function (sometimes called [faj]).

\subsection{} We are interested in this note in the case $k=1$, for which of course one cannot apply the above formulas. The first one ``suggests'' that the main term vanishes, and that we should expect an asymptotic of the size of $\sqrt{q}$. More precisely, the error term in (\ref{1}) comes from a sum involving class numbers of quadratic fields: we will be surprised to see that this fact is preserved! One purpose of this note is to make these heuristics a bit more precise.

\subsection{Construction of dihedral weight one modular forms}\label{induction} For simplicity, let's assume the level $q$ is prime in this section only, and let $\chi_q=\left(\frac{\cdot}{q}\right)$ be the Legendre symbol (hence $q\equiv 3\,(4)$ is the parity condition stated above). One can construct cusp forms of arbitrary weight using induction of Hecke characters associated to quadratic fields; see Miyake \cite{M1} \S 4.8. Let $K_q=\Q(\sqrt{-q})$, and let $\psi$ be a character of the class group of $K$. Let:
\begin{eqnarray}
\theta_\psi(z):=\sum_{\n\subset \entier}\psi(\n)\e(N(\n)z)
\end{eqnarray}
Then $\theta_\psi$ is actually a weight one primitive cusp form with nebentypus $\chi_q$, unless $\psi$ is real (in which case, the theta series is a weight one Eisenstein series). They are paired two by two, and one gets $\frac{h(K)-1}{2}$ independent such forms. The construction can be carried out for ray class groups, and for real quadratic fields as well (but in the last case, one gets a non-squarefree level).\\

Let's call $\mathcal{S}_1^{Dih}(q,\chi)$ the span of these theta series; by Siegel's theorem (for the lower bound):
$$ q^{1/2-\varepsilon} \ll_\varepsilon{\rm dim}\bigg(\mathcal{S}_1^{Dih}(q,\chi)\bigg)=\frac{h(K)-1}{2}\ll_\varepsilon q^{1/2+\varepsilon}$$
which gives the $\sqrt{q}$ term promised earlier. The big task is to understand what are the forms we missed in this construction, and why they should be rare (though it is still conjectured). For the moment, let me state the first form of:

\begin{conj}\label{conj} One has the estimate, for $q$ varying among the squarefree integers:
$${\rm dim}\bigg(\mathcal{S}_1(q,\chi)\bigg)=\frac{h(K_q)}{2}+\mathcal{O}_\varepsilon(q^\varepsilon)$$
with an $\mathcal{O}$-constant independent of $q,\chi$.
\end{conj}

\subsection{Remark.} From the point of view of representation theory, the theta series is really ``induction'' in the usual sense. Indeed, $\psi$ defines a character of the Galois group $G_K=\Gal{K}$ by class field theory, hence a one-dimensional representation. So its induction ${\rm Ind}_{G_K}^{G_\Q}\psi$ defines a two-dimensional complex Galois representation. It is conjectured by Langlands that such operations on the Galois side still correspond to automorphic forms: in this case, the theta series \emph{is} the corresponding form. For more on this issue, the paper of Rogawski \cite{Rog} is excellent; Gelbart \cite{G} sketches the proof in general for quadratic induction; Bump \cite{Bu} has a very clear introduction to Langlands conjectures; Miyake \cite{M1} gives a proof in the classical setting; etc.

\subsection{How to tackle the conjecture?} By using the trace formula, the best estimate of the dimension of weight one modular forms is about $q$, perhaps with a power saving in $\log(q)$. Duke achieved a power saving in $q$, by using two conflicting properties:
\begin{itemize}
\item ``Smoothness'' via harmonic analysis: a general feature of automorphic forms, a bit degenerated in weight one.
\item ``Rigidity'' via Galois representations, imposing strong limitations.
\end{itemize}

The bound Duke gets is roughly of  size  $q^{11/12}$, using the full strength of Serre's analysis of the possible lifts from $\rm PGL(2)$ to $\rm GL(2)$; if one is happy with a tiny power saving, one can soften  the input a lot.

\section{Harmonic analysis}

I will be brief here: it is difficult to give a better survey than Michel's Park City lecture notes \cite{M}. Plus, I need only the large sieve inequality -- which is the small residue of harmonic analysis left in weight one. The book of Iwaniec and Kowalski \cite{IK} contains many examples of such inequalities as well, and explains how one can guess the ``best'' bound.

\subsection{} Let $\mathcal{B}_1(q,\chi)$ be an \emph{orthonormal} basis of $\mathcal{S}_1(q,\chi)$. Recall that $c_n(f)$ denotes the $n$-th Fourier coefficient of such a form. 
\begin{proposition}
Let $\{x_n\}_{n\in\N}$ be an arbitrary sequence of complex numbers. The following estimate is uniform in $N, q, \chi$:
\begin{eqnarray}\label{sieve}
\sum_{f\in\mathcal{B}_1(q,\chi)}\bigg|\sum_{n=1}^N x_nc_n(f)\bigg|^2\ll\bigg(1+\frac{N}{q}\bigg)\sum_{n=1}^N|x_n|^2
\end{eqnarray}
\end{proposition}
The estimate (\ref{sieve}) is called \emph{large sieve inequality}: it is a powerful substitute for Cauchy-Schwarz inequality. Three proofs are given in Michel's lectures: using Kuznetsov trace formula, Rankin-Selberg $L$-functions, and a variant of Iwaniec (explained in Duke's paper as well). The first one usually gives the best bounds, but requires that the family one averages over be spectrally complete. Instead of reproducing these proofs here, let's explain  this spectral aspect.

\subsection{}\label{kuznetsov} The proof of (\ref{sieve}) should start as follows: one expands the product on the left hand side; one needs then to estimate the average of the product of two Fourier coefficients over an orthonormal basis. Such a formula is a ``Petersson-Kuznetsov''-type trace formula. For example, for weights greater than 2, one has (the implied constant depends on $k$ and can be made explicit) the Petersson formula:
$$\sum_{f\in\mathcal{B}_k(q,\chi)}c_n(f)\overline{c_m(f)}\doteq \delta_{n,m}+\sum_{c\equiv 0\,(q)}\frac{K\ell_\chi(n,m,q)}{|c|}J_{k-1}\bigg(\frac{4\pi\sqrt{nm}}{|c|}\bigg)$$
with a Bessel function $J_{k-1}$ and a twisted Kloosterman sum $K\ell_\chi$, which is an explicit exponential sum, and can be evaluated for example using Weil's bound.\\

 Now let $\mathcal{M}_k(q,\chi)$ be an orthonormal basis of the space of cuspidal Maass cusp forms of weight $k$, nebentypus $\chi$ , and $\Phi$ a ``nice'' holomorphic, rapidly decreasing test function, one has the Kuznetsov formula:

$$\sum_{f\in\mathcal{M}_k(q,\chi)}\Phi(t_f)c_n(f)\overline{c_m(f)}+({\textrm{Eis.}})\doteq \delta_{n,m}+\sum_{c\equiv 0\,(q)}\frac{K\ell_\chi(n,m,q)}{|c|}\widehat{\Phi}\bigg(\frac{4\pi\sqrt{nm}}{|c|}\bigg)$$
where: $1/4+t_f^2$ is the Laplace eigenvalue of $f$, $\widehat{\Phi}$ is a Bessel-transform of $
\Phi$, (Eis.) denotes the contribution of the continuous spectrum parametrized by the Eisenstein series (for the definition fo Maass forms, their multipliers, etc., see Bump \cite{Bu} or Michel \cite{M}). To ensure the space of weight $k$ Maass forms is not trivial, there is a parity-type condition on $k,\chi$. Interestingly, holomorphic forms can be embedded into the space of Maass form, using:
$$f\in\mathcal{S}_k(q,\chi)\longmapsto y^{k/2}f(z).$$
It turns out that for $k\geq 2$, the image of this map is precisely the space of forms with Laplace eigenvalue $\frac{k}{2}\big(1-\frac{k}{2}\big)$, which is an isolated point of the spectrum of the Laplacian. By choosing a relevant test function (i.e. taking $\Phi$ with support outside the vertical line parametrized by the continuous spectrum), one can recover Petersson formula. \\

For $k=1$, $1/4$ is no longer away from the continuous spectrum, but one can still prove a ``large sieve inequality'' for weight one holomorphic forms using such an embedding, and a positivity argument (a sum  of positive numbers  on $\mathcal{B}_1(q,\chi)$ can certainly be bounded by such a sum on $\mathcal{M}_1(q,\chi)!$). The first impression is that such a naive argument loses a lot of information. However, by re-establishing good harmonic behaviour, one can actually improve Duke's bound, without changing the input: I will come back to this later.

\subsection{Remark.} If $\chi$ is an even ($\chi(-1)=1$) character, the non-holomorphic Maass cusp forms with eigenvalue $1/4$ are conjectured to correspond to Galois representations, as do their ``odd'' cousins, the holomorphic weight one modular forms.

\subsection{Rankin-Selberg convolution}\label{rankin} Let $f,g$ be two modular forms with the same weight (the construction works for any weight, but requires a slightly different normalization), possibly with distinct levels and nebentypus. We saw how to define a reasonable theory of $L$-functions, so let's suppose that $f,g$ are both primitive newforms (or at least Hecke eigenforms at the good primes). One then defines (absolute convergence for $\Re(s)\gg 0$):
$$L(s,f\times\bar{g})=\zeta(\chi_f\overline{\chi_g},2s)\sum_{n\geq 1}\frac{a_n(f)\overline{a_n(g)}}{n^s}$$
This $L$-function satisfies a functional equation, has a simple pole at $s=k$ if and only if $f=g$ (up to a scalar multiple), and has an Euler product, whose local factors at good primes are given by:
$$L(s,f_p\times \overline{g_p})=\prod_{(i,j)\in\{1,-1\}}(1-\alpha_i(p,f)\overline{\alpha_j(p,g)}p^{-s})^{-1}.$$
As $L(s,f\times\bar{f})$ has non-negative coefficients, Landau's lemma implies the convergence of the series in the region $\{\Re(s)>k\}$, from which one can infer that $|\alpha_i(p,f)|< p^{k/2}$ (by the absence of poles in this region applied to the factor at $p$) -- and many more important results, like $\displaystyle{\sum_{n=1}^N |c_n(f)|^2\ll_{\varepsilon,f} N^{k+\varepsilon}}$ for any $\varepsilon>0$, which is Ramanujan on average, cf. Michel on this issue.\\

In the region of absolute convergence (and ignoring the bad primes, finite in number) $$\displaystyle{\log\big(L(s,f\times\bar{f})\big)=\sum_{m=1}^\infty g_m(s)}$$ where
$$g_m(s)=-\sum_{(p,q)=1}\sum_{(i,j)\in\{1,-1\}}\frac{\alpha_i(p,f)^m\overline{\alpha_j(p,f)}^mp^{-ms}}{m}$$
one has
$$\sum_{m=1}^\infty g_m(s)=-\log(s-k)+\mathcal{O}_{s\to k}(1)$$
and so
$$g_1(s)=\sum_{(p,q)=1}|a_p(f)|^2p^{-s}<-\log(s-k) {\textrm{ for }s\to k}$$

\subsection{Remark.} Once the Ramanujan conjecture is proven, it is easy to see that $\sum_{n\geq 2}g_m(s)$ stays bounded at $1$, so the above inequality is actually an asymptotic.

\subsection{Consequence.}\label{cons} Suppose $f$ is a weight one modular form and an eigenvector of the good Hecke operators. Then, the Hecke eigenvalues are integers of a finite extension $\Q(f)$ of $\Q$, so for any real number $M$ the set:
$$S(M)=\{\alpha\in \mathcal{O}_{\Q(f)}\,:\, |\sigma(\alpha)|\leq M\, \forall \sigma: \Q(f)\hookrightarrow \C\}$$
is finite. So $\{a_f(p)\,:\, p\not| q\ \textrm{ and } a_f(p)\in S(M)\}$ is finite. Also one has for $\Re(s)>1$:
$$\sum_{p\not\,| q}\frac{|a_f(p)|^2}{p^s}\geq \sum_{\substack{p\not\,| q\\ a_f(p)\notin S(M) }} \frac{M^2}{p^{s}}$$
and so:
$$\sum_{\substack{p\not\,| q\\ a_f(p)\notin S(M) }}p^{-s} < -M^{-2}\log(s-1)$$
which, by the definition of the (upper) Dirichlet density, means:
$${\textrm {Dens}}(p\not| q\,:\, a_f(p)\notin S(M))\leq M^{-2}$$
and we can summarize:
\begin{lemme}
For any positive $\eta$, there exists a finite subset $Y_\eta$ of $\C$, and a set of primes $\mathscr{P}_\eta$ of density less than $\eta$ such that:
$$\forall p\notin\mathscr{P}_\eta\,,\,a_f(p)\in Y_\eta. $$
\end{lemme}  
This means that, possibly outside a small set of primes, the Hecke eigenvalues of weight one modular forms (with fixed weight, nebentypus) are finite in  number. Deligne-Serre  rules out this possible tiny set.

\subsection{Densities.} A subset $\mathscr{P}$ of prime numbers has upper Dirichlet density $\alpha$ if:
$$\limsup_{s\to 1}\frac{\sum_{p\in\mathscr{P}}p^{-s}}{\sum_{p\textrm{ prime}}p^{-s}}=\alpha.$$
If the above quotient actually has a limit, one talks about ``Dirichlet density''. The (upper) natural density would be defined by:
$$\limsup_{x\to\infty}\frac{1}{\pi(x)}\{p\in\mathscr{P}, p\leq x\}$$
again with the symbol ``$\lim$'' when it makes sense. Its existence implies the existence of Dirichlet density. For integers, one has of course the two corresponding definitions, again with the same comparisons. We will always write ${\rm Dens}(\mathscr{P})$ to denote the (upper) Dirichlet or natural density; as soon as the Tchebotarev density theorem is involved, natural densities actually exist.

\section{Galois representations.}\label{galois}
Conseils de lecture: I already mentionned  Rogawski's review \cite{Rog} on Artin conjecture , Langlands-Tunnell theorem. I should add Taylor's article \cite{T}, very complete but more difficult. Bump \cite{Bu} has in his first chapter an excellent introduction to the topic, and it is not possible to me to forget Rohrlich's paper \cite{Roh} on the Weil-Deligne group and elliptic curves. Many other things are to be found in Cornell-Silverman-Stevens \cite{CS}. Bushnell and Henniart's book \cite{BH} on the local Langlands conjectures is also a good companion.\\

Let $k$ be a complete field (topologized with an absolute value, allowed to be trivial only if $k$ is a finite field), $V$ a finite dimensional topological vector space on $k$ (so topologically $V\cong k^n$ for some $n$). Let $K$ be a field. A representation of the absolute Galois group $\Gal{K}=G_K$ in $V$ is a pair $(\rho,V)$, often denoted $\rho$, where $V$ is as above, and $\rho$ is a continuous homomorphism:
$$\rho: G_K\longrightarrow {\rm GL}(V).$$
$\bullet\,\, \rho$ is \emph{irreducible} if there are no closed subspace of $V$ invariant under $\rho$.\\
$\bullet\,\, \rho$ is \emph{semi-simple} if $V$ is a direct sum of irreducible closed subspaces.\\
We will be interested in the cases where $k$ is $\C$, a $\p$-adic field or a finite field, with its natural topology.

\subsection{Remark.} Choose an algebraic closure $\bar{k}$ of $k$, and suppose that there are countably many finite extensions of $k$ (inside this choice of $\bar{k}$). Then I claim that there exists a finite extension $E$ of $k$ such that $\rho(G_K)\subset \GL_n(E)$ -- i.e. $\rho$ is rational. Indeed, $\rho(G_K)$ is compact, and one can write
$$\rho(G_K)=\bigcup_{[E:k]<\infty}\bigg(\rho(G_K)\cap \GL_n(E)\bigg) .$$
By the Baire category theorem ($E$, as a finite dimensional $k$-vector space, is endowed with the product topology hence is complete), one of the $\rho(G_K)\cap \GL_n(E)$ has non-empty interior in $\rho(G_K)$, so is cofinite inside it; by adding the remaining matrices, one gets the claim.

\subsection{The semi-simplicity issue.}\label{semisimplification} $G_K$ is always compact. So if $k=\C$, any representation is semi-simple (one can define the topology of $V$ by a $G_K$-invariant product, and can thus take orthogonal complements; this fact extends of course in the case of Hilbert spaces).\\
In the other cases, it is a bit more complicated. But if $G$ is a finite group, then any representation of $G$ is semisimple, as long as ${\rm char}(k)$ is prime to $|G|$; in particular it is always true for $\p$-adic fields.\\
The importance of semi-simple representations is that they are determined (up to isomorphism) by their character, or characteristic polynomial:
\begin{proposition}\label{ss}
Let $G$ be a group, $\rho,\rho'$ two $k$-semi-simple representations. If for any $g\in G$, ${\rm Tr}(\rho)(g)={\rm Tr}(\rho')(g)$ for ${\rm char}(k)=0$ (resp. ${\rm det}(Id-X\rho(g))={\rm det}(Id-X\rho'(g))$ if $k$ is finite), then $\rho,\rho'$ are isomorphic.
\end{proposition}
In the finite case, the proof can be found in Curtis-Reiner \cite{CR} (30.16) e.g. Look in Husem\"oller's elliptic curves \cite{H} \S 15.2 for some comments in the $\ell$-adic case.\\

If $\rho$ is not semi-simple, then one can define from $\rho$ a semi-simple representation (unique up to isomorphism). Indeed, let $\{V_i\}$ be a finite Jordan-H\"older composition series of $V$: $\{0\}=V_0\subset V_1\subset \dots \subset V_N=V$ where all inclusions are strict, and $V_i$ is a maximal $G$-submodule of $V_{i+1}$. Consider ${\rm Gr}(V)=\bigoplus_{i\geq 0}V_{i+1}/V_i$, which carries a $G$-action, induced from the one on $G$: this representation is semi-simple by construction, and has the same character, characteristic polynomials, as $\rho$. It is called \emph{the} semi-simplification of $\rho$ (by the ``uniqueness'' implied by the last proposition).

\subsection{Complex case.} $G_\Q$ has a basis of neighbourhoods of the identity consisting of subgroups, namely $\{{ Ga\ell}(\overline{\Q}/E)\}_{[E:\Q]<\infty}$. But $\GL_n(\C)$, as any Lie group, has a neighbourhood of the identity which doesn't contain any subgroup: this is easily seen by looking at ${\rm det}: \GL_n(\C)\rightarrow \C^\times$, as if $U\subset \GL_n(\C)$ is small enough, ${\rm det}(U)$ is contained in a small disc centered at 1, which cannot contain any subgroup of $\C^\times$. If $U$ is chosen this way, then for $E/\Q$ large enough, $\rho({ Ga\ell}(\overline{\Q}/E))\subset U$, so is trivial! As a consequence, any complex Galois representation factors through some finite Galois extension $E$ of $\Q$, which can of course be chosen so that $\rho: {Ga\ell}(E/\Q)\rightarrow \GL_n(\C)$ is injective.

\subsection{If $K=\Q$.} Let $(\rho,V)$ be a complex or $\p$-adic Galois representation. Let us fix once and for all an algebraic closure of $\Q, \Q_p$ for all the primes $p$, and an embedding $\overline{\Q}\hookrightarrow \overline{\Q_p}$ for each $p$. From these embeddings, one deduces a system of local representations:
$$\rho_p:G_{\Q_p}\longrightarrow \GL(V).$$
Recall the short exact sequence
$$\{1\}\longrightarrow I_p\longrightarrow G_{\Q_p}\longrightarrow G_{\mathbf{F}_p}=\widehat{\Z}\longrightarrow \{0\}$$
where $\widehat{\Z}$ is topologically generated by the Frobenius $x\mapsto x^p$. One can lift this Frobenius to an element of $G_{\Q_p}$, defined up to the inertia subgroup $I_p$. Let's fix a choice for each $p$ of such a Frobenius element noted $F_p$. It defines an operator on inertia-invariant vectors $V^{I_p}$, and so one can define, if $k=\C$
$$L(s,\rho_p):={\rm det}\bigg(Id-X\rho(F_p)\big| V^{I_p}\bigg)^{-1}{}_{\big| _{X=p^{-s}}}$$
It is important to remark that this $L$-factor does not depend on the choice of $\overline{\Q}\hookrightarrow \overline{\Q_p}$.

If $k$ is $\p$-adic, one has to choose first an algebraic embedding $\imath:k\hookrightarrow \C$, and then take
$$L(s,\rho_p):=\imath\bigg({\rm det}\big(Id-X\rho(F_p)\big| V^{I_p}\big)^{-1}\bigg){}_{\big| _{X=p^{-s}}}$$
Note that the last $L$-factor depends on the algebraic embedding $\imath$, but fortunately Deligne tells us that in geometric situations, the relevant determinant lies in a number field.\\
The global $L$-function, called the Artin $L$-function, is defined for $\Re(s)\gg 0$ by
$$L(s,\rho)=\prod_pL(s,\rho_p)$$
Actually, if $k=\C$, the eigenvalues of the Frobenius elements have modulus one, so it is easy to check the covergence for $\Re(s)>1$; in the $\p$-adic case, I don't know in general: if $\rho$ arises in the \'etale cohomology of an algebraic variety, then Deligne's theorem (Weil I or II) gives a bound on the eigenvalues of the geometric Frobenius elements, from which one can deduce the convergence.

\subsection{Remark.} In the complex case, $\rho$ factors through a finite quotient, which already kills almost all inertia subgroups; one deduces that almost all the local $L$-factors have degree $n$ in $p^{-s}$ -- one says that $\rho$ is \emph{unramified} at $p$. In the $\p$-adic case, it can be false in general; but if $\rho$ is geometric, then Grothendieck \emph{et al.} showed that $\rho$ actually factors through a quotient ${Ga\ell}(E/\Q)$, usually infinite, such that $E/\Q$ is unramified outside a finite number of primes. This is nice, because one can use a powerful tool: the Tchebotarev theorem.

\subsection{Tchebotarev density theorem} It makes possible a  refinement of proposition \ref{ss}. Let's state it in the most general form:
\begin{theoreme}
Let $E/\Q$ be an algebraic Galois extension, unramified outside a finite set of primes. Let $C$ be a measurable subset of $Ga\ell(E/\Q)$ stable under conjugation, with ${\rm vol}(\partial C)=0$. Then:
$${\rm Dens}\bigg(p{\textrm{ unramified prime }}\,:\, F_p\subset C\bigg)=\frac{{\rm vol}(C)}{{\rm vol}(Ga\ell(E/\Q))}$$
\end{theoreme}

Here ``density'' means natural or Dirichlet density; note that in this statement, one considers the Frobenius elements as conjugacy classes, as one did not fix a priori any embedding $\overline{\Q}\hookrightarrow \overline{\Q_p}$: as a consequence, one has to take into account the action of $Ga\ell(E/\Q)$ on the decomposition group (by conjugation; a transitive action), and the Frobenius only defines a conjugacy class. We'll talk about``finite'' (resp. ``infinite'') Tchebotarev if $E/\Q$ is a finite (resp. infinite) extension. \\

Let $g\in Ga\ell(E/\Q)$. For any finite Galois subextension $K/\Q$, one can apply Tchebotarev to the conjugacy class of the reduction of $g$ modulo $Ga\ell(K/\Q)$: it tells that this reduction $g_K$ belongs to the conjugacy class $F_p$ for infinitely many primes $p$. But one the other hand, $g_K\to g$\footnote[2]{to make  this convergence meaningful, consider the filter of the sections for the order defined by the inclusion of subextensions.}, so \emph{the Frobenius classes are dense in $Ga\ell(E/\Q)$}. Another way to deduce this from the above version of Tchebotarev, is to take $C=\{hUh^{-1}\,:\,h\in G_{\Q}\}$ for an  open neighbourhood of $g$ of the type $U=g\cdot Ga\ell(E/K)$. In any case, one deduces the sought amelioration of prop. \ref{ss}:
\begin{proposition}\label{cheb}
Suppose $k$ is $\C$, a finite or a $\p$-adic field. Let $\rho,\rho'$ be two $k$-semi-simple representations of $G_\Q$. If for any unramified $p$, oustide a zero-density set of primes, ${\rm Tr}\big(\rho(F_p)\big)={\rm Tr}\big(\rho'(F_p)\big)$ for ${\rm char}(k)=0$ (resp. ${\rm det}\bigg(Id-X\rho(F_p)\bigg)={\rm det}\bigg(Id-X\rho'(F_p)\bigg)$ if $k$ is finite), then $\rho,\rho'$ are isomorphic.
\end{proposition}
The case where $k$ is a finite field case is not completely justified, but here again the Galois representation factors through a finite quotient, and one needs only the ``finite'' Tchebotarev theorem.

\subsection{Remark.} It is not easy to find a reference with a proof of the above version, which I stated for purely aesthetic reasons, as the condition on the boundary of $C$ is quite restrictive in practice; see for example Serre \cite{S2}. One can actually deduce it from the ``finite'' version, amply covered in the literature. Indeed, the ``finite'' version implies that for any subset $C$ stable by conjugation, 
$${\rm Dens}\bigg(p{\textrm{ unramified prime }}\,:\, F_p\subset CGa\ell(E/K)\bigg)=\frac{|\pi_K(C)|}{|Ga\ell(K/\Q)|}$$
for any finite (Galois) subextension $K/\Q$, where we noted $\pi_K:Ga\ell(E/\Q)\twoheadrightarrow Ga\ell(K/\Q)$ the projection. Then, as the Haar measure of the (infinite) Galois group is the projective limit of the measures of its finite quotients, one sees easily that
$${\rm vol}(\overline{C})=\lim_{K}\frac{|\pi_K(C)|}{|Ga\ell(K/\Q)|}$$
hence the condition on the boundary of $C$. But to actually justify that the primes in the theorem have a density, it is necessary to argue as follows. Let $Ga\ell(E/\Q)^{\#}$ be the set of  conjugacy classes of $Ga\ell(E/\Q)$. Let $\pi:Ga\ell(E/\Q) \twoheadrightarrow Ga\ell(E/\Q)^\#$ the (open) projection. The Haar measure defines a measure on the conjugacy classes (by pushforward); the ``finite'' Tchebotarev, together with the above fact on projective limits of measures, implies that for any $K/\Q$ finite, and any $c\in Ga\ell(E/\Q)^\#$
$$\frac{\log(X)}{X}\sum_{p\leq X}\delta_{F_p}({\rm Char}_{cGa\ell(E/K)})_{\substack{\longrightarrow \\ x\to\infty}}\int_{Ga\ell(E/\Q)^\#}{\rm Char_{cGa\ell(E/K)}}(g)d^\#g$$
But as the characteristic functions ${\rm Char}_{cGa\ell(E/K)}$ are dense (in the space of continuous functions on $Ga\ell(E/\Q)^\#$), the above convergence is true for any continuous map (uniform boundedness theorem), hence the sequence of the Frobenius elements $\{F_p\}_{p}$ is equidistributed in $Ga\ell(E/K)^\#$. From this, as one can approximate the characterisitic function of any measurable set $C$ in $Ga\ell(E/\Q)^\#$ -- with negligible boundary -- by continuous functions (here the approximation is in the sense of the order on the space of real functions, not in norm, cf. Billingsley \cite{B} theorem 2.1, or Serre \cite{S2} appendix of chapter I), one gets the result. \\
The point is that the densities (natural or Dirichlet) on prime numbers are not a measure, and the naive argument (``take the intersection of the sets of primes'') doesn't work; but one gets the result by equidistribution, in two steps: from special characteristic functions to continuous functions (by approximation in norm), then from continuous functions to more general characteristic functions (thanks to the order on real numbers).

\subsection{} If $\rho:G_\Q\longrightarrow GL(V)$ is a complex or $\p$-adic semi-simple Galois representation, we defined above its Dirichlet $L$-series by a convergent (in some half plane) Euler product. It turns out that this $L$-function satisfies a functional equation, of the type:
$$\Lambda(s,\rho)=q^{\frac{1-2s}{2}}\omega_\rho\Lambda(1-s,\rho^{\vee })$$
where $\Lambda$ is the product of $L$ with some $\Gamma$-factors, describing the behaviour of $\rho$ at infinity (cf. the Durham volume, Martinet (same volume as \cite{S1}), and Tate in the Corvallis volume \cite{Ta}); $\rho^\vee$ is the contragredient of $\rho$, $\omega_\rho$ is a complex number of modulus one. Brauer showed that this $L$-function extends to a meromorphic function on $\C$, and Artin conjectured that $\Lambda$ is holomorphic unless $\rho$ contains the trivial representation: cf. Heilbronn in Cassels-Fr\"ohlich \cite{CF} for the meromorphic continuation.

\subsection{Remark.} There is at least one other type of geometric $L$-function one should mention: the $L$-function obtained from algebraic varieties by ``counting points'' modulo a prime. For example, for abelian varieties, it turns out that this zeta function agrees with the $L$-function of the Galois action on the (dual of the) Tate module: cf. Rohrlich \cite{Roh} for a proof of this in the case of elliptic curves, and much more.

\subsection{Remark.} Unlike complex Galois representations, usually $\p$-adic Galois representations don't have finite image.

\section{The Deligne-Serre theorem}

\subsection{Statement of the theorem.}
\begin{theoreme}\label{deligne-serre}
Let $f$ be a weight one cuspidal modular form, with nebentypus $\chi$, and level $q$, a Hecke eigenvector at the good primes. Then there exists an irreducible complex Galois representation, unique up to isomorphism, $\rho_f:G_\Q\longrightarrow \GL(2,\C)$ such that at any good prime $p$:
$$L(s,\rho_p)=L(s,f_p).$$
\end{theoreme}

\begin{corollaire}[The Ramanujan Conjecture]\label{ramanujan}
Under the same hypotheses, at good primes one has:
$$|a_p(f)|\leq 2$$
If furthermore $f$ is a newform, one has:
$$|a_n(f)|\leq \tau(n)\ll_\varepsilon n^\varepsilon.$$
If $f$ is any weight one modular form, then:
$$|c_n(f)|\ll_{\varepsilon,f}n^\varepsilon.$$
\end{corollaire}

\begin{corollaire}\label{artin}
Under the same hypotheses, if $f$ is a primitive newform, the local $L$ factors actually agree on the ramified primes as well, the conductors of the Artin $L$-function and of the newform $f$ agree, as do the $\varepsilon$ factors.
\end{corollaire}

\noindent {\sc Proof} of Corollary \ref{ramanujan}: Firstly $a_p(f)={\rm Tr}(\rho_f(F_p))$; as $\rho_f(G_\Q)$ is finite, the eigenvalues of $\rho_f(F_p))$ are roots of unity, hence the first point; for the second, note that by multiplicativity $a_n(f)=\prod_pa_{p^{v_p(n)}}(f)$, and finally $a_{p^k}(f)=\prod_{i+j=k}\alpha_1(p,f)^i\alpha_2(p,f)^j$ whence the second bound. The third is trivial from the decomposition fo the total space of modular forms into newforms, and oldforms. {\sc q.e.d.}\\

See Deligne-Serre for the proof of the second claim: by the theorem, the two $L$-functions agree except possibly on finitely many ramified primes; so the quotient of the $L$-functions is a finite product. Using the functional equation, one proves it easily by reductio ad absurdum.\\

The proof of the Ramanujan conjecture is much more difficult in weight $k\geq 2$, as the geometric nature of the Galois representation attached to $f$ is needed to prove it, via the study of the weights in $\ell$-adic cohomology (Deligne, Weil I).

\subsection{Langlands' conjectures} These predict that Galois representations should ``correspond'' to certain adelic automorphic forms. The Deligne-Serre theorem belongs to this program. Recently, Khare \emph{et al.} proved that reciprocally, any odd Galois representation is associated to a weight one modular form (odd means that the Dirichlet character $\chi$, corresponding to ${\rm det}(\rho(.))$ by class field theory, satisfies $\chi(-1)=-1$), completing the Langlands-Tunnell theorem (cf. Rogawski \cite{Rog}). With a little more work, one could strengthen the theorem with Langlands' philosophy in mind. If $f$ is a  newform, one can lift $f$ to a square-integrable (modulo the center) function on $\GL(2,\Q)\backslash\GL(2,\adele_\Q)$; let $V(f)=\overline{{\rm span}\{f(\cdot\, g)\,:\, g\in\GL(2,\adele_\Q)\}}$: it turns out that $V(f)$ is an irreducible unitary representation $\pi_f$ of $\GL(2,\adele_\Q)$ -- cf. Gelbart for a justification of these statements. Such a representation splits as a completed tensor product of  unitary representations of the $p$-adic Lie groups $\GL(2,\Q_p)$:

$$\pi_f\cong \widehat{\bigotimes_{p\leq \infty}}\pi_{f,p}.$$
The local Langlands conjectures require that the $\pi_{f,p}$'s correspond to representations $\rho_p $ of the Weil-Deligne group of $\Q_p$. The local/global compatibility we had in mind is, with $\mathcal{W}_{\Q_p}\subset G_{\Q_p}$ the absolute Weil group:
$${\rho_f}_{\big|_{\mathcal{W}_{\Q_p}}}={\rho_p}_{\big|_{\mathcal{W}_{\Q_p}}}.$$
%which makes sense because the corresponding Weil-Deligne representations in this case are of Galois type.

\subsection{Consequence.}\label{linear} Let $f, \rho_f$ be as in the theorem; $\rho_f(G_\Q)$ is a finite subgroup of $\GL(2,\C)$: so is $\overline{\rho_f(G_\Q)}$ in ${\rm PGL}(2,\C)$, which is thus conjugate to one of the following: $\Z/n\Z, D_{2n}, \mathfrak{A}_5, \mathfrak{A}_4, \mathfrak{S}_4$. The first is ruled out, because it's commutative ($\rho_f$ would be a sum of two characters in this case); it is an interesting fact, not too difficult to establish, that the dihedral groups are obtained precisely from modular forms induced from characters on quadratic extensions, as explained in section \ref{induction}. One calls a form $f$ \emph{exotic} if $\overline{\rho_f(G_\Q)}$ is isomorphic to one of $\mathfrak{A}_5, \mathfrak{A}_4, \mathfrak{S}_4$. For $\mathcal{B}^{new}_1(q,\chi)$ a basis of $\mathcal{S}^{new}_1(q,\chi)$, with forms chosen to be Hecke eigenforms, let
$$s^{Exotic}(q,\chi)=|\{f\in\mathcal{B}^{new}_1(q,\chi)\,:\, f \textrm{ is exotic }\}|.$$
This number is easily seen to be independent of the choice of $\mathcal{B}^{new}_1(q,\chi)$. 
\begin{conj}\label{conjecture} With the notations above:
$$s^{Exotic}(q,\chi)\ll_\varepsilon q^\varepsilon$$
\end{conj}

As the conjecture, and its weakened forms, are statements uniform in $q$, it is important to see what happens when $q$ varies. For example, if it is true that for $q$ fixed, the number of Hecke eigenvalues of dihedral forms is bounded (they are sums of two roots of unity of order determined by $q$), when $q$ varies, the number of possibilities grows with $q$. But, one can expect that for exotic forms, one can get a uniform bound on the number of Hecke eigenvalues.\\

Indeed, suppose that $(q,\chi)$ varies among pairs of integers and Dirichlet characters having order \emph{bounded} by $M$. Let $f$ be an exotic primitive (new)form. The image of $\rho_f(G_\Q)$ in ${\rm PGL}(2,\C)$ has order bounded by 120; so the eigenvalues of $\rho(F_p)$ modulo the center are roots of unity of order less than 120, so the eigenvalues of $\rho(F_p)$ itself are roots of unity of order less than $120M$; so the $a_f(p)$, when $p$ varies among the primes, is a sum of two such complex numbers. Let $S$ denote  this finite set. As a consequence, the polynomial 
$P(X)=\prod_{s\in S}(X-s)+A$
satisfies, for any prime $p$, any level $q$, and any exotic form $f$:
$$P(a_p(f))=A$$
Using the Hecke-type relation at good primes:
$$a_p(f)^2=a_{p^2}(f)+\chi(p)$$
expanding the product giving $P(a_p(f))$, and choosing $A$ large enough, one gets a \emph{linear} relation of the type:
$$\sum_{k=1}^{|S|}\mu_k(\chi)a_{p^k}(f)=B(\chi)$$
valid for any $q,\chi$, and any good prime $p$, with bounded (when $\chi,q$ vary) complex numbers $\mu_k(\chi)$, $|B(\chi)|\geq 1$. This rough linear relation is sufficient to get a power saving as a first attempt to prove the (still unsolved) conjecture. The point is that one can suppress the assumption on the boundedness of the order of $\chi$ to get a linear relation among the $a_{p^k}(f)$, and even get an explicit (small) length, but this requires a careful and difficult analysis of the possible lifts from $\rm PGL(2)$ to $\GL(2)$. This linear relation depends on the type of the form: in the case of icosahedral forms $f$ (i.e. $\overline{\rho_f(G_\Q)}\simeq \mathfrak{A}_5$)  one has for example:
$$\bar{\chi}^6(p)a_{p^{12}}(f)-\bar{\chi}^4(p)a_{p^8}(f)-\bar{\chi}(p)a_{p^2}(f)=1$$
The other relations are similar, but involve smaller powers of $p$, and so the bound obtained is better -- cf. section \ref{duke}.

\section{Proof of the Deligne-Serre theorem.}
The reader should refer to the original paper, and Serre's overview in the Durham volume. We will content ourselves to recall the main steps. The starting point is the following theorem (due to Deligne \emph{et al.}):

\begin{theoreme}\label{deligne}
Let $f$ be a weight $k\geq 2$ cuspidal modular form, with nebentypus $\chi$, level $q$, Hecke eigenvector at the good primes with eigenvalues $\{a_p(f)\}_{(p,q)=1}$. Let $K$ be a number field containing $\Q(f)$, and let $\lambda$ be a finite prime of $K$ of residue characteristic $\ell$.
Then there exists an irreducible $\lambda$-adic Galois representation, unique up to isomorphism, $\rho_f:G_\Q\longrightarrow \GL(2,K_\lambda)$, unramified away from $q\ell,$ such that at any prime $p$ \emph{not dividing $q\ell$:}
\begin{eqnarray}
{\rm Tr}(\rho_f(F_p))=a_f(p)\\
{\rm det}(\rho_f(F_p))=p^{k-1}\chi(p)
\end{eqnarray}
\end{theoreme}

Once Deligne-Serre  is proven, the above  is true for weight one as well, by choosing an algebraic embedding $\imath: K_\lambda\hookrightarrow \C$. The uniqueness statement follows from Tchebotarev as we saw in section \ref{galois} -- this is the reason why one does not state the compatibility using the $L$-factors, which depend on the algebraic embedding, but this could be fixed by requiring that $\imath_{|_{\Q(f)}}={\rm id}_{\Q(f)}$.

\subsection{Step 1: Reduction modulo a prime.} Let's consider the following preliminary result:

\begin{theoreme}\label{reduction}
Let $f$ be a weight $k\geq 1$ cuspidal modular form, with nebentypus $\chi$, level $q$, Hecke eigenvector at the good primes with eigenvalues $\{a_p(f)\}_{(p,q)=1}$. Let $K$ be a number field containing $\Q(f)$, and let $\lambda$ be a finite prime of $K$ of residue characteristic $\ell$ and let $k_f={\mathbf{F}_\ell}(a_f(p), \chi(p); p\not |q)$.
Then there exists a semi-simple Galois representation, unique up to isomorphism, $\overline{\rho_f}:G_\Q\longrightarrow \GL(2,k_f)$, unramified away from $q\ell,$ such that at any prime $p$ \emph{not dividing $q\ell$:}
\begin{eqnarray}
{\rm Tr}(\overline{\rho_f}(F_p))\equiv a_f(p)\, {\rm mod}(\lambda)\\
{\rm det}(\overline{\rho_f}(F_p))\equiv p^{k-1}\chi(p)\,{\rm mod}(\lambda).
\end{eqnarray}
\end{theoreme}

One can strengthen the theorem, as it is actually only necessary that $f$ be a Hecke eigenform modulo $\lambda$: see Deligne-Serre. 
\begin{enumerate}
\item If $f$ is a weight one form, multiplication of $f$ by a weight $m$ $\rm SL(2,\Z)$ Eisentein series produces a weight $m+1$ modular form. By considering the normalized Eisenstein series:
$$E_m(z)=1-\frac{b_m}{2m}\sum_{n=1}^\infty \sigma_{m-1}(n)\e(nz),$$
the Bernouilli numbers $b_m$ satisfy congruence relations (Clausen-von Staudt theorem): $\ell b_m\equiv -1\,{\rm mod}(\ell)$ if $(\ell-1)|m$ (cf. Borevich-Shafarevitch, chap. 5 \S 8), so $fE_m\equiv f \,{\rm mod}(\lambda)$ for such a $m$.
\item Re-establish good Hecke behaviour: there exists a weight $m+1$ modular form $f'$, eigenform at good primes, defined over some extension $K'/K$, such that $a_p(f')\equiv a_p(f)\,{\rm mod}(\lambda')$ for some $\lambda'|\lambda$, and $p\not| \ell q$.
\item One can apply Deligne's theorem \ref{deligne} to this $f'$, getting a representation $\rho'$ with values in $\GL(2,K'_{\lambda'})$.
\item $\rho'(G_\Q)$ is a compact subgroup of $\GL(2,K'_{\lambda'})$, which one can suppose to be contained in $\GL(2,{\mathcal{O}'}_{\lambda'})$, after possible conjugation.  So one can reduce $\rho'$ modulo $\lambda'$, and get $\overline{\rho'}$  with values in $\GL(2,\mathcal{O}'_{\lambda'}/\lambda')$. Note that so far we preserved the congruences between the initial $f$ and the Galois representations, 
\item but $\overline{\rho'}$ may not be semi-simple: taking its semi-simplification (cf. \ref{semisimplification}), one gets a semi-simple representation satisfying the congruences required in the theorem.
\item One can reduce the field of definition from $\mathcal{O}'_{\lambda'}/\lambda'$ to $k_f$, because the characteristic polynomials of all the $\overline{\rho'}(g)$ have their coefficients in $k_f$ by definition of $k_f$ (note that Tchebotarev theorem has been used again): more precisely, $\overline{\rho'}$ is isomorphic to a representation $\overline{\rho_f}:G_\Q\longrightarrow \GL(2,k_f)\subset \GL(2,\mathcal{O}'_{\lambda'}/\lambda')$ (cf. Weintraub, Theorem 6.5, \S 3.6).
\end{enumerate}

The unramifiedness is of course evident. One has lost  irreducibility of the representations by reducing them. The study of (odd) Galois representations with values in finite fields led Serre to conjecture that they all arise as described above, which was proven recently by Khare \emph{et al.}.

\subsection{Step 2: Bounding the images of the reductions.}
In step one, one reduced Galois representations modulo a prime. To get a complex Galois representation, one needs to lift them back! This process is not trivial, but one can do it in the following situation: 
\begin{proposition}\label{lift}
Let $G$ be a finite group, $\ell$ be a prime not dividing $|G|$ and $\overline{\rho_\ell}: G\longrightarrow \GL(n,\F_\ell)$ a representation. Suppose there exists a number field $K$ containing the roots of unity of order $|G|$, and a prime ideal $\lambda$ for which $\entier/\lambda=\F_\ell$. Then there exists a representation $\rho: G\longrightarrow \GL(n,\entier)$ whose reduction modulo $\lambda$ is isomorphic to $\overline{\rho_\ell}$ .
\end{proposition}

Cf. Weintraub, end of \S 3.6, for a proof.\\

So, if $K\supset \Q(f)$ is ``sufficiently large'', one decomposes $\overline{\rho_f}$ of  theorem \ref{reduction} into $G_\Q\longrightarrow \overline{\rho_f}(G_\Q)$ followed by $\overline{\rho_f}(G_\Q)\hookrightarrow \GL(2,k_f)$, and would like to lift the second to $\GL(2,K)$. The result we quoted tells us this is possible, if we can assume that the characteristic of $k_f$ is coprime to $|\overline{\rho_f}(G_\Q)|$, at least for sufficiently many primes $\lambda$ of $K$. This is indeed the case, and one proceeds as follows.

Let $K\supset \Q(f)$ be a finite Galois extension of $\Q$; let $\mathscr{S}$ be the set of primes of $\Q$, totally split in $K$. It is known that ${\rm Dens}(\mathscr{S})=1/[K:\Q]$, so 
$\mathscr{S}$ is infinite, and we have some flexibility in the choice of $K$. For $\ell\in\mathscr{S}$, theorem \ref{reduction} gives a semi-simple representation $\overline{\rho_\ell}: G_\Q\longrightarrow \GL(2,\F_\ell)$. As a consequence, $G_\ell:=\overline{\rho_\ell}(G_\Q)$ is a semi-simple subgroup of $\GL(2,\F_\ell)$.\\

The key result is the following:
\begin{proposition}
With this notation, one has:
$$\sup_{\ell\in\mathscr{S}}|G_\ell|<\infty$$
\end{proposition}

\noindent{\sc Sketch of proof:} Recall that in section \ref{cons}, we proved the existence of a set of primes $\mathscr{P}_\eta$, for any positive $\eta$, having density less than $\eta$, and a finite set of complex numbers $Y_\eta$ such that:
$$\forall p\not\in\mathscr{P}_\eta, a_p(f)\in Y_\eta.$$
Let $H_{\ell,\eta}$ be the set containing all the conjacy classes $\{\overline{\rho_\ell}(F_p)\}_{p\not\in\mathscr{P}_\eta}$. By Tchebotarev, one has:
$${\rm Dens}(\mathscr{P}^c_\eta)=\frac{|H_{\ell,\eta}|}{|G_\ell|}\geq 1-\eta.$$
In other words, $H_{\ell,\eta}$ is a subset ``large'' inside $G_\ell$, whose set of characteristic polynomials is small (less than $|Y_\eta|$, but in any case independent of $\ell$). As the semi-simple subgroups of ${\GL}(2,\F_\ell)$ can be classified, one checks case by case the uniform bound claimed in the proposition: cf. Deligne-Serre, proposition 7.2, for the details.

\subsection{Step 3: Lift to $\GL(2,\C)$.}
One is now in position to apply proposition \ref{lift}. Indeed, let $M=\sup_{\ell\in\mathscr{S}}|G_\ell|$: by adjoining them if necessary, on can suppose that $K$ contains the $M$-th roots of unity; and let 
$$\mathscr{S}'=\{p\textrm{ prime }: p \textrm{ totally split in K and } p> M\}.$$
Proposition \ref{lift} produces, for each $\ell\in\mathscr{S}'$, a representation $\rho_\ell: G_\Q\longrightarrow \GL(2,\entier)$ whose reduction modulo a prime $\lambda$ dividing $\ell$ is $\overline{\rho_\ell} $ ($\ell$ being totally split in $K$, any prime ideal $\lambda$ dividing $\ell$ satisfies $\entier/\lambda\cong \F_\ell$).\\

The last point to be careful with is the effect of this lifting on the congruences as stated in theorem \ref{reduction}. One can argue as follows. 
Consider the finite set of polynomials:
$$\mathcal{P_M}=\{P(X)=(1-\alpha X)(1-\beta X):\, \alpha,\beta\textrm{ $M$-th roots of unity}\}$$
\begin{enumerate}
\item For each prime $p$, and $\ell\in\mathscr{S}', \ell\not=p$ and $\lambda|\ell$, the theorem \ref{reduction} implies that:
$$\exists R_\ell\in\mathcal{P_M}\,:\, 1-a_p(f)X+\chi(p)X^2\equiv R_\ell(X)\,{\rm mod}(\lambda).  $$
As $\mathcal{P_M}$ is finite, the same polynomial works for infinitely many $\ell$, so $1-a_p(f)X+X^2$ itself is in $\mathcal{P_M}$. Note that at this point, the Ramanujan conjecture is proven.
\item As the set $\mathcal{P}_M$ is finite, one can suppose that $\forall P,Q\in\mathcal{P}_M; P\not=Q\Rightarrow P\not\equiv Q\,{\rm mod}(\lambda)$ (after removing a finite set of unsatisfying primes $\ell$).
\item Fix a prime $\ell$ not dividing the level $q$ of the form $f$, and let $p\neq\ell$ as above. The lift $\rho_\ell$ to $\GL(2,\entier)$ is unramified outside $q\ell$, and the characteristic polynomial of $\rho_\ell(F_p)$ is in $\mathcal{P_M}$: this is because by construction $|\rho_\ell(G_\Q)|=|G_\ell|$, so $\rho_\ell(F_p)$ has order less than $M$. On the other hand, the polynomial $1-a_p(f)X+\chi(p)X^2$ is in $\mathcal{P_M}$ as well, and as $\rho_\ell$ lifts $\overline{\rho_\ell}$, the two polynomials are congruent modulo
 $\lambda$: by the choice of $\ell$, they are hence equal, so one has:
 $$\forall p\not|q\ell, {\rm det}(Id-X\rho_\ell(F_p))=1+a_p(f)X+\chi(p)X^2.$$
 \item One has the same conclusion for another choice $\ell'$ (by replacing $\ell$ by $\ell'$).
 \item This means that the two representations $\rho_\ell, \rho_{\ell'}:G_\Q\longrightarrow \GL(2,\C)$ are isomorphic as $\C$-representations, because of proposition \ref{cheb} (after choosing an embedding $K\hookrightarrow \C$), that actually all the $\rho_\ell$'s are unramified outside $q$, and that the relation:
  $$\forall p\not|q, {\rm det}(Id-X\rho_\ell(F_p))=1+a_p(f)X+\chi(p)X^2$$
holds. After this choice, one renames $\rho_\ell$ into $\rho_{f}$, and the construction is finished.
\end{enumerate}

One should note that the finiteness of $\rho_f(G_\Q)$ is a consequence of the construction; and that one gets the continuity of $\rho_f$ from the finiteness of the image!

\subsection{Step 4: Irreducibility.} Reductio ad absurdum: one doesn't construct a stable subspace, but infers that reducibility would contradict the bound given on the Rankin-Selberg convolution seen in \ref{rankin}. Indeed, if $\rho_f$ were not irreducible, by semi-simplicity (automatic over $\C$), one could find two one-dimensional stable subspaces of $\C^2$, on which $G_\Q$ would act by characters $\chi_1,\chi_2$, which by class field theory one can view as Dirichlet characters:
$$\rho_f\cong \chi_1\oplus\chi_2$$
So $a_p(f)=\chi_1(p)+\chi_2(p)$, and $\chi(p)={\rm det}(\rho_f(F_p))=\chi_1(p)\chi_2(p)$, which proves  by the way that $\chi_1\not=\chi_2$ (else $\chi(-1)=1$ and $f=0$). One gets
$$\sum_{p\not\,| q}\frac{|a_p(f)|^2}{p^s}=2\sum_{p\not\,| q}p^{-s}+\sum_{p\not\,| q}\frac{\chi_1(p)\overline{\chi_2(p)}}{p^s}+\sum_{p\not\,| q}\frac{\chi_2(p)\overline{\chi_1(p)}}{p^s}$$
It is well known that for any non-trivial Dirichlet character $\psi$, $\sum_n \psi(n)n^{-s}$ is holomorphic at $s=1$, so the two last terms are bounded. This would imply that
$$\sum_{p\not\,| q}\frac{|a_p(f)|^2}{p^s}=2\sum_{p\not\,| q}p^{-s}=-2\log(s-1)+\mathcal{O}_{s\to 1}(1).$$
This contradicts the result given in section \ref{rankin}, and proves therefore the irreducibility of $\rho_f$.

\section{The dimension of the space of weight one modular forms.}\label{duke}
Let's go back to the statements we gave on the dimension of weight one forms in the conjectures \ref{conjecture} (or \ref{conj}, given in the case of squarefree levels, in which case the main term is explicitly related to quadratic inductions). We wrote
$${\rm dim}(\mathcal{S}^{new}_1(q,\chi))=s^{Dih}(q,\chi)+s^{Exotic}(q,\chi)$$
and expect that for any $\varepsilon>0$
$$s^{Exotic}(q,\chi)\ll_\varepsilon q^\varepsilon$$
Duke \emph{et al.} \cite{Du} proved the following:
\begin{proposition}\label{duke1}
There exists a positive real number $\delta$ such that:
$$s^{Exotic}(q,\chi)\ll q^{1-\delta}.$$
Explicitly, the above inequality holds for any $0\leq\delta < \frac{1}{12}$.
\end{proposition}

\noindent{\sc Proof:} 
\subsection{} The idea is the following: let $\{x_n\}_{n\geq 1}$ be a sequence of complex numbers. By positivity one has the following inequality:
$$\sum_{f\in\mathcal{B}^{new}_1(q,\chi)}\bigg|\sum_{n=1}^N x_nc_n(f)\bigg|^2\leq \sum_{f\in\mathcal{B}_1(q,\chi)}\bigg|\sum_{n=1}^N x_nc_n(f)\bigg|^2$$
for any orthonormal basis $\mathcal{B}^{new}_1(q,\chi)$  of the space of newforms, completed into an orthonormal basis $\mathcal{B}_1(q,\chi)$ of the total space space of weight one forms $\mathcal{S}_1(q,\chi)$. To use Galois representations, one chooses $\mathcal{B}^{new}_1(q,\chi)$ to be made of Hecke eigenforms, normalized ($a_1(f)=1$) so that the Hecke eigenvalues and Fourier coefficients agree. This normalization destroys the unitarity of the forms, so if $\mathcal{B}^{new}_1(q,\chi)$ denotes now a Hecke eigenbasis of primitive forms, one has to study
$$\sum_{f\in\mathcal{B}^{new}_1(q,\chi)}\frac{1}{\|f\|^2}\bigg|\sum_{n=1}^N x_na_n(f)\bigg|^2.$$
\subsection{} One makes three packets in $\mathcal{B}^{new}_1(q,\chi)$: octahedral, icosahedral, tetrahedral according to the type of the projective representation $\overline{\rho_f}$ ($\mathfrak{S}_4, \mathfrak{A}_5, \mathfrak{A}_4$ respectively). For each type, one chooses the test vector by using a linear relation, seen in section \ref{linear}, of the type
$$\sum_{k=1}^{|S|}\mu_k(\chi)a_{p^k}(f)=B(\chi)$$
The longest linear relation is obtained for icosahedral forms, namely:
\begin{eqnarray}\label{lin}\bar{\chi}^6(p)a_{p^{12}}(f)-\bar{\chi}^4(p)a_{p^8}(f)-\bar{\chi}(p)a_{p^2}(f)=1
\end{eqnarray}
for which one chooses
\begin{displaymath}
x_n = \left\{ \begin{array}{ll}\bar{\chi}^6(p) & \textrm{ if }n=p^{12}, (p,q)=1\\-\bar{\chi}^4(p) & \textrm{ if }n=p^8, (p,q)=1\\-\bar{\chi}(p) & \textrm{ if }n=p^2, (p,q)=1\\ 0 & \textrm{ otherwise}  \end{array}\right.
\end{displaymath}
For this choice of $\{x_n\}$, one has for any icosahedral forms
$$ \bigg|\sum_{n=1}^N x_na_n(f)\bigg|^2 =\big|\{p\textrm{ prime };\, (p,q)=1\textrm{ and } p^{12}\leq N \}\big|^2\gg \left(\frac{N^{\frac{1}{12}}}{\log(q)}\right)^2$$
by the prime number theorem; so one gets
$$\sum_{f\in\mathcal{B}^{ico}_1(q,\chi)}\frac{1}{\|f\|^2} \left(\frac{N^{\frac{1}{12}}}{\log(q)}\right)^2\ll\sum_{f\in\mathcal{B}^{ico}_1(q,\chi)}\frac{1}{\|f\|^2}\bigg|\sum_{n=1}^N x_na_n(f)\bigg|^2$$
\subsection{} On the other hand, the right hand side is bounded above by the large sieve inequality, which gives here
$$\sum_{f\in\mathcal{B}^{ico}_1(q,\chi)}\frac{1}{\|f\|^2}\bigg|\sum_{n=1}^N x_na_n(f)\bigg|^2\ll \bigg(1+\frac{N}{q}\bigg)\left(\frac{N^{\frac{1}{12}}}{\log(q)}\right)$$
\subsection{}\label{a} The two last steps give
$$ \sum_{f\in\mathcal{B}^{ico}_1(q,\chi)}\frac{1}{\|f\|^2} \ll \bigg(1+\frac{N}{q}\bigg)\left(\frac{N^{\frac{1}{12}}}{\log(q)}\right)^{-1}$$
which is uniform in $q,\chi,N$; one still has to choose $N$, which depends on the behaviour of $\|f\|$ for a primitive form of level $q$.
\subsection{} Such an estimate is classical: cf. \cite{I} section 13.6, and \cite{M} section 2.3. Indeed:
\begin{multline}\nonumber\|f\|^2:=\int_{\Gamma_0(q)\backslash \mathfrak{H}}|\sqrt{y}f(z)|^2 \frac{dxdy}{y^2}\\\doteq {\rm res}_{s=1}{\rm vol}(\Gamma_0(q)\backslash \mathfrak{H})\int_{\Gamma_0(q)\backslash \mathfrak{H}}|\sqrt{y}f(z)|^2E_s(z)\frac{dxdy}{y^2}
\end{multline}
 where the implied constant in ``$\doteq$'' is absolute, and with the usual Eisenstein series 
$$E_s(z)=\sum_{\gamma\in\Gamma_\infty\backslash\Gamma_0(q)}\frac{y^s}{|cz+d|^{2s}}$$
By unfolding the Eisenstein series, one proves that:
$$ \int_{\Gamma_0(q)\backslash \mathfrak{H}}|\sqrt{y}f(z)|^2E_s(z)\frac{dxdy}{y^2}=L(s,f\times\bar{f})$$
On the other hand, the Rankin-Selberg $L$-series (or its adelic counterpart) factorizes as\footnote[2]{For the definition of the adjoint $L$-function  $L(s,f,{\rm Ad})$, see \cite{Rog} section 13 e.g.}
$$L(s,f\times\bar{f})\doteq L(s,1)L(s,f,{\rm Ad})$$
which shows that
$${\rm res}_{s=1}L(s,f\times\bar{f})\doteq L(1,f,{\rm Ad})$$
and, at this point, one uses the well-known bound (cf. \cite{M}, section 1.3)
$$L(s,f,{\rm Ad})\ll_\varepsilon q^\varepsilon$$
This gives:
$$\|f\|^2\ll_\varepsilon {\rm vol}(\Gamma_0(q)\backslash \mathfrak{H})q^\varepsilon\ll_\varepsilon q^{1+\varepsilon}$$
\subsection{} By plugging the last estimate in \ref{a}, one has
$$|\mathcal{B}_1^{ico}(q,\chi)|\ll_\varepsilon q^{1+\varepsilon} \bigg(1+\frac{N}{q}\bigg)\left(\frac{N^{\frac{1}{12}}}{\log(q)}\right)^{-1}$$
\subsection{} By choosing $N=q$, one gets the result for icosahedral forms, i.e.:
$$|\mathcal{B}_1^{ico}(q,\chi)|\ll_\varepsilon q^{1-\frac{1}{12}+\varepsilon}$$
One does the same for the other types of forms -- the only difference is that the linear relation among the Hecke eigenvalues have shorter length, and give a better estimate (better towards the conjecture).

\subsection{Remark.} If one uses the very weak linear relation obtained in section \ref{linear}, one gets the estimate -- cf. \ref{linear} for the notations:
$$s^{Exotic}(q,\chi)\ll_\varepsilon q^{1-|S|^{-1}+\varepsilon}$$
which is still a power saving, and uses much less analysis as the one required to get the explicit identity (\ref{lin}).

\subsection{Remark.} It is striking that one can ameliorate the bound of proposition \ref{duke1}, without changing the input: Michel-Venkatesh \cite{MV} noted that the use of Kuznetsov trace formula, instead of the large sieve inequality, gives a bound in $q^{6/7}$ -- instead of $q^{11/12}$ here. Here is how they proceed: embed the space of weight one modular forms into the space of weight one Maass cuspforms, as explained in section \ref{kuznetsov}, and complete the basis $\mathcal{B}_1^{new}(q,\chi)$ into a basis $\mathcal{M}_1^{new}(q,\chi)$ of Maass cuspforms. One gets, by choosing a test function $\Phi$ non-negative on the spectrum of the hyperbolic Laplacian:
$$\sum_{f\in\mathcal{B}^{new}_1(q,\chi)}\Phi(0)\bigg|\sum_{n=1}^N x_nc_n(f)\bigg|^2\leq \sum_{f\in\mathcal{M}_1(q,\chi)}\Phi(t_f)\bigg|\sum_{n=1}^N x_nc_n(f)\bigg|^2$$
One still has the lower bound:
$$\sum_{f\in\mathcal{B}^{new}_1(q,\chi)}\bigg|\sum_{n=1}^N x_nc_n(f)\bigg|^2\gg_\varepsilon \left(\frac{N^{\frac{1}{12}}}{\log(q)}\right)^2 q^{-1-\varepsilon}s^{Exotic}(q,\chi)$$
To get an upper bound, one expands the square, and uses the Kuznetsov formula, to get
\begin{multline}\nonumber
\sum_{f\in\mathcal{M}_1(q,\chi)}\Phi(t_f)\bigg|\sum_{n=1}^N x_nc_n(f)\bigg|^2\doteq\sum_{n=1}^N|x_n|^2\\+\sum_{n,m=1}^Nx_n\overline{x_m}\sum_{c\equiv 0\,(q)}\frac{K\ell_\chi(n,m,q)}{|c|}\widehat{\Phi}\bigg(\frac{4\pi\sqrt{nm}}{|c|}\bigg)\\ \ll_{\varepsilon,\Phi} \sum_{n=1}^N|x_n|^2+(qN)^\varepsilon\frac{\sqrt{N}}{q}\left(\sum_{n=1}^N|x_n|\right)^2
\end{multline}
where the last estimate follows from Weil's bound on Kloosterman sums. By choosing the same test vector $\{x_n\}$ as before, one gets the bound
$$s^{Exotic}(q,\chi)\ll_\varepsilon q^{\frac{6}{7}+\varepsilon}$$

\section{On the Fourier coefficients of weight one modular forms.}
Our last application of Deligne-Serre theorem is an estimate of the number of non-vanishing Fourier coefficients of weight one modular forms. Once one can use the associated Galois representation, the argument is incredibly simple and short, and is a nice interaction between worlds of Galois representations and of the analytic theory of modular forms. The problem is the following:
\begin{probleme}
Let $f\in\mathcal{S}_1^{new}(q,\chi)$. What is the behaviour, when $N\to\infty$, of the quantity 
$\frac{1}{N}\bigg|\bigg\{n\leq N\,:\, a_n(f)\not=0\bigg\}\bigg|$?
\end{probleme}
The idea is the following: one first studies  the question over prime integers, for which the number $a_p(f)$ is related to the trace of a matrix; then one bootstraps the estimate to get the density over the integers. The last section of Deligne-Serre article discusses the problem. Let $f\in\mathcal{S}_1^{new}(q,\chi)$, $\rho_f$ its associated Galois representation given by theorem \ref{deligne-serre}, and let $E_f/\Q$ be the finite Galois extension whose Galois group is $G_\Q/{\rm ker}(\rho_f)$ .

\subsection{} For any prime $p$, one has $a_p(f)={\rm Tr}(\rho_f(F_p))$. Let:
$$Z(\rho_f)=\{g\in Ga\ell(E_f/\Q)\,:\, {\rm Tr}(\rho_f(g))=0\}$$
One can apply Tchebotarev's theorem, to get:
$${\rm Dens}\bigg(\Big\{p\textrm{ prime }\,:\, a_p(f)=0\Big\}\bigg)={\rm vol}(Z(\rho_f))=\beta>0$$
The number $\beta$ is non-zero, because the (class of the) complex conjugation belongs to $Z(\rho_f)$.
\subsection{} One uses Wirsing's theorem:
\begin{theoreme}
Let $\omega$ be a non-negative multiplicative function, such that:
\begin{eqnarray}
\sum_{p\leq N}\omega(p)\log(p)= \Big(\alpha+o_{\scriptscriptstyle{N\to\infty}}(1)\Big)N\\
\exists \gamma_1,\gamma_2>0\,:\,\forall p \textrm{ prime and }k\geq 0,  f(p^k)\leq \gamma_1\gamma_2^k
\end{eqnarray}
Then one has:
$$\sum_{n\leq N}\omega(n)=\Big(1+o_{\scriptscriptstyle{N\to\infty}}(1)\Big)\cdot\frac{N}{\log(N)}\cdot\frac{\exp(-\gamma\alpha)}{\Gamma(\alpha)}\cdot\prod_{p\leq N}\bigg(\sum_{k=0}^\infty \frac{\omega(p^k)}{p^k}\bigg)$$
where $\gamma$ denotes the Euler constant.
\end{theoreme}
If one applies the theorem for $\omega(n)=1$ if $a_f(n)\not= 0$, $\omega(n)=0$ otherwise, one gets:
$$\bigg|\Big\{n\leq N\,:\, a_n(f)\not=0\Big\}\bigg|=\Big(C+o_{\scriptscriptstyle{N\to\infty}}(1)\Big)\times \frac{N}{\log(N)^\beta}$$
for $C$ some real number (depending on $f$), as the product over the primes has an asymtotic contribution for the primes such that $a_f(p)\not=0$, which is a set of density $\alpha=1-\beta$.

\subsection{Remark} One can rephrase these results as follows: for  a fixed primitive eigenform of weight one $f$, $h:\R \rightarrow \C$ a  bounded continuous function:
$$\frac{1}{N}\sum_{n=1}^N h(a_n(f))_{\substack{ \longrightarrow \\ {N\to\infty}}} h(0)$$
For a fixed form $f$, the set of possible values of the $a_p(f)$ is finite: let $\{a_i\}_i$ denote this set, and let $\beta_i$ be the density of primes $p$ for which $a_p(f)=a_i$. Then, one has for any continuous function $h:[-2,2]\rightarrow \C$:
\begin{eqnarray}\label{sato1}
\frac{\log(N)}{N}\sum_{p=2}^N h(a_p(f))_{\substack{ \longrightarrow \\ {N\to\infty}}} \sum_{i}\alpha_i h(a_i) 
\end{eqnarray}
which shows that the sequence of the probability measures associated to the Hecke eigenvalues weakly converges to the discrete measure appearing on the right hand side. This situation is diametricly opposed to what is expected to hold for higher weight modular forms:

\subsection{The Sato-Tate conjecture.}
\begin{conj}[Sato-Tate]
Let $f$ be a (non-{\sc cm}) primitive newform of weight $k \geq 2$, level $q$. Let $\{a_p(f)\}_{p\textrm{ prime}}$ denote its Hecke eigenvalues. For any continuous function $h:[-2,2]\rightarrow \C$, one has:
\begin{eqnarray}\label{satok}
\frac{\log(N)}{N}\sum_{p=2}^N h\left(\frac{a_p(f)}{p^{\frac{k-1}{2}}}\right) {}_{\substack{ \longrightarrow \\ {N\to\infty}}}\frac{1}{\pi}\int_{-2}^2 h(u)\sqrt{1-\frac{u^2}{4}}du.
\end{eqnarray}
\end{conj}
Cf. Ribet \cite{R1} for a definition of non-{\sc cm} in this setting. The difference between the limit distributions appearing in (\ref{sato1}) (discrete atomic probability), and (\ref{satok}) (given by a density), is another instance of the dichotomy weight one/weight greater than one. Another difference is that Sato-Tate conjecture has been proven (recently)  for some non-{\sc cm} forms of weight two, corresponding to elliptic curves without complex multiplication. If one takes Sato-Tate conjecture for granted, one deduces easily that:

$${\rm Dens}\bigg(\Big\{p\textrm{ prime }\,:\, a_p(f)=0\Big\}\bigg)=0$$
I would like to sketch here an unconditional proof of the above equality, as it illustrates the power of Galois representations -- when one has them at one's disposal!\\

So let $f$ be a (cuspidal) primitive form of weight $k\geq 2$ and level $q$. If $K$ is a number field containing $\Q(f)$, and $\lambda|\ell$ a prime of $K$ of residue characteristic $\ell$ not dividing the level, theorem \ref{deligne} gives a $\lambda$-adic representation:
$$\rho_{f,\lambda}:G_\Q \longrightarrow \GL(2,K_\lambda)$$
 One can assume that $\rho_{f,\lambda}(G_\Q)\subset \GL(2,\mathcal{O}_\lambda)$ (after conjugating if necessary). Let $\overline{\rho_{f,\lambda}}$ be the (semi-simplification of the) reduction modulo $\lambda$ of $\rho_{f,\lambda}$. Then, trivially, one has:
$$Z(f):=\Big\{p\textrm{ prime}:\, a_p(f)=0\Big\}\subset Z(f,\lambda):= \Big\{p\textrm{ prime}:\, a_p(f)\equiv 0\,{\rm mod}(\lambda)\Big\}$$
Let $E_\lambda/\Q$ be the finite Galois extension corresponding to the sugbroup ${\rm ker}(\overline{\rho_{f,\lambda}})\subset G_\Q$; by Tchebotarev, one has:
$${\rm Dens}\bigg(Z(f,\lambda)\bigg)=\frac{|\{ g\in \overline{\rho_{f,\lambda}}(G_\Q)\,:\, {\rm Tr}(g)\equiv 0 \}|}{|\overline{\rho_{f,\lambda}}(G_\Q)|}$$
Suppose for the moment that one knows that $\rho_{f,\lambda}(G_\Q)= \GL(2,\mathcal{O}_\lambda)$ (which happens for almost all primes in the case of a non-{\sc cm} elliptic curve). Then, by letting $q_\ell=|\mathcal{O}_\lambda/\lambda|$, one has:
$$|\overline{\rho_{f,\lambda}}(G_\Q)|=(q_\ell^2-1)(q_\ell^2-q_\ell)\,;\,|\{ g\in \overline{\rho_{f,\lambda}}(G_\Q)\,:\, {\rm Tr}(g)\equiv 0 \}|=(q_\ell^2-1)(q_\ell-1) $$
which gives:
$$\lim_{\ell\to\infty}{\rm Dens}\bigg(Z(f,\lambda)\bigg)=0$$
Unfortunately, the images of Galois representations corresponding to modular forms are not always open. Ribet studied in different works this problem, and gave in \cite{Ribet} a general result, valid for non-{\sc cm} modular forms. For our purposes, it is sufficient to say that for \emph{totally split} primes $\ell$, $$\rho_{f,\lambda}(G_\Q)=\{g\in\GL(2,\Z_\ell)\,:\, {\rm det}(g)\in \Z_\ell^{\times(k-1)}\}$$ Here $\Z_\ell^{\times(k-1)}$ is the sugbroup of $k-1$ powers in $\Z_\ell^{\times}$. 
The subgroup $\overline{\rho_{f,\lambda}}(G_\Q)$ is open and cofinite inside $\GL(2,\Z_\ell)$, with index independent of $\ell$,  and so one can argue as we did above! Indeed, $[\Z_\ell^\times:\Z_\ell^{\times( k-1)}]=\frac{|\mu_{k-1}(\Q_p)|}{|k-1|_p}$ stays bounded when $\ell\to\infty$ among totally split primes; so $[\GL(2,\Z_\ell):\rho_{f,\lambda}(G_\Q)]$ does as well (the two last indexes are equal to each other); as a consequence $[\GL(2,\F_\ell):\overline{\rho_{f,\lambda}(G_\Q)}]$ ($\leq [\GL(2,\Z_\ell):\rho_{f,\lambda}(G_\Q)]$) is uniformly bounded as the totally split primes $\ell$ go to infinity, and we have:
$$\frac{|\{ g\in \overline{\rho_{f,\lambda}}(G_\Q)\,:\, {\rm Tr}(g)\equiv 0 \}|}{|\overline{\rho_{f,\lambda}}(G_\Q)|} \ll_{k,f} \frac{(\ell^2-1)(\ell-1)}{(\ell^2-1)(\ell^2-\ell)}$$

This shows unconditionally, that:
$${\rm Dens}\bigg(\Big\{p\textrm{ prime }\,;\, a_p(f)=0\Big\}\bigg)=0$$
in the case of a non-{\sc cm} modular form of weight $k\geq 2$.

\subsection{Remark} We have not dealt with forms having complex multiplication, but one expects a  limit density; for {\sc cm} elliptic curves, the limit distribution is the sum of a  measure with  density on $[-2,2]$ (pushforward of the uniform measure on the circle) and a Dirac at zero (coming from inert primes), and this is unconditionally proven (cf. Murty in \cite{CKM}, Lecture 1, end of \S 2).

\subsection{A conjecture} Let's suppose that $q$ is a prime number here, and that $\chi_q$ is  Legendre character. Let $f$ be a primitive Hecke eigenform in $\mathcal{S}_1(q,\chi_q)$. We proved in (\ref{sato1}) the existence of numbers $\{\alpha_i(f)\}_{1\leq i\leq N_f}$ with sum 1, and of real numbers $\{a_i(f)\}_{1\leq i\leq N_f}$ in $[-2,2]$ such that:
$$\frac{1}{\pi(N)}\sum_{p=2}^N \varphi(a_p(f)){}_{\substack{ \longrightarrow \\ {N\to\infty}}} \mu_f\Big(\varphi\Big):= \sum_{i=1}^{N_f}\alpha_i(f) \varphi\big(a_i(f)\big)$$
for any  continuous function $\varphi:[-2,2]\rightarrow \C$. The limit distribution $\mu_f$ is discrete, but one can wonder if one re-establishes a density after averaging over a basis of newforms. 

\begin{conj}\label{conjecture2}
Let $\mathcal{B}_1(q)$ be the basis of primitive newforms of $\mathcal{S}_1(q,\chi_q)$. Then the sequence of probability measures $$ \mu_q:=\frac{1}{{\rm dim}\bigg(\mathcal{S}_1(q,\chi_q)\bigg)}\sum_{f\in\mathcal{B}_1(q) } \mu_f$$ converges weakly as the prime $q\equiv 3\, {\rm mod}(4)$ goes to infinity. More precisely, for any continuous function $\varphi:[-2,2]\rightarrow \C$, one has:
$$ \frac{1}{{\rm dim}\bigg(\mathcal{S}_1(q,\chi_q)\bigg)}\sum_{f\in\mathcal{B}_1(q) } \mu_f\Big(\varphi\Big)\,{}_{\substack{ \longrightarrow \\ {q\to\infty}}} \,\frac{\varphi(0)}{2}+\frac{1}{4\pi}\int_{-2}^2 \varphi(x)\frac{dx}{\sqrt{1-\frac{x^2}{4}}} $$
\end{conj}

If conjecture \ref{conjecture} is true, or if  $s_1^{Exotic}(q,\chi_q)=o_{\scriptscriptstyle{q\to\infty}}\Big(s^{Dih}_1(q,\chi_q)\Big)$, then conjecture \ref{conjecture2} is plausible, as the only contribution comes from dihedral forms. Assuming that the torsion in the class groups of imaginary quadratic fields is ``small'' (cf. \cite{EV} for the terminology), we can then show that the mean value of such class group characters are equidistributed in the unit circle, as the discriminant tends to infinity. The Dirac comes from inert primes: if $f$ is dihedral, then the primes inert in $\Q(\sqrt{-q})$ form a set of density $1/2$, and they all have $a_p(f)=0$.\\

An even more challenging question is the following (``horizontal'' Sato-Tate). One fixes a prime $p$,  and considers for a continuous function $\varphi$ as above the expression
$$\mu_{p,q}(\varphi):=\frac{1}{{\rm dim}\bigg(\mathcal{S}_1(q,\chi_q)\bigg)}\sum_{f\in\mathcal{B}_1(q) } \varphi(a_p(f))$$
Serre studied the weak limit of this sequence of measures, when the weight is greater than two, as the level $q$ tends to infinity, using Eichler-Selberg's trace formula (see section 29 of \cite{KL}). It is plausible that in the harmonically degenerated case of weight one forms, there is no limit anymore; but it would be interesting to determine  the weak cluster points, and it is possible that the smoothed sequence $\widetilde{\mu}_{p,Q}:=\frac{\log(Q)}{Q}\sum_{q\leq Q}\mu_{p,q}$ has a limit when the prime $Q$ tends to infinity. \\

Finally, the depth aspect. Fix an imaginary quadratic field and consider a family of ray class characters whose conductors are unbounded, or similarly a family of ring-class characters associated to a decreasing family of orders: in any of these examples, when the characters induce a weight one modular form of a certain level, it is possible to compare the size of the dihedral dimension to the exotic dimension, and to look at the weak limits of the corresponding families of measures.

Stanford University, Department of Mathematics, building 380, Stanford, California 94305, USA.\\
E-mail: \href{mailto:trotabas@math.stanford.edu}{\texttt{trotabas@math.stanford.edu}}

\end{document}